\documentclass[11pt]{article}
\usepackage[utf8]{inputenc}
\usepackage[english]{babel}
\usepackage{amsmath}
\usepackage{amsfonts}
\usepackage{amssymb}
\usepackage{graphicx}
\usepackage{fancyhdr}
\usepackage{enumitem}
\usepackage{listings}
\usepackage{booktabs}
\usepackage{color}
\usepackage{moresize} 
\usepackage{caption}
\usepackage{epstopdf}
\usepackage{bm}
\newtheorem{theorem}{\bf Theorem}[section]

\newtheorem{es}{\bf Example}
\newcommand{\Proofend}{\hfill$\diamondsuit$}
\def\PP{{\mathbb P}}
\def\C{{\mathcal{C}}}

\newcommand{\al}{{\alpha}}

\def\NN{{\mathbb N}}

\newcommand{\wal}{w_\alpha}

\newcommand{\uu}{u}

\title{De la Vall\'ee Poussin filtered polynomial approximation on the half line}
\author{D.Occorsio, W.Themistoclakis}

\begin{document}

\maketitle
\begin{abstract}
On the half line we introduce a new sequence of near--best uniform approximation polynomials, easily computable by the values of the approximated function at a truncated number of Laguerre zeros. Such approximation polynomials come from a discretization of filtered Fourier--Laguerre partial sums, which are filtered  by using a de la Vall\'ee Poussin (VP) filter. They have the peculiarity of depending on two parameters: a truncation parameter  that determines how many of the $n$ Laguerre zeros are considered, and a localization parameter, which determines the range of action of the VP filter that we are going to apply. As $n\to\infty$, under simple assumptions on such parameters and on the Laguerre exponents of the involved weights, we prove that the new VP filtered approximation polynomials
have uniformly bounded Lebesgue constants and uniformly convergence at a near--best approximation rate, for any locally continuous function on the semiaxis.
\newline
The theoretical results have been validated by the numerical experiments. In particular, they show a better performance of the proposed VP filtered approximation versus the truncated Lagrange interpolation at the same nodes, especially  for functions a.e. very smooth with isolated singularities. In such cases we see a more localized approximation as well as a good reduction of the Gibbs phenomenon.
\end{abstract}

\section{Introduction}
The paper deals with the approximation of  locally continuous functions $f$ defined on the positive half-line. The function $f$ we consider can  exponentially grow at infinity or/and have an algebraic singularity at the origin. Measuring the approximation error with a suitable weighted uniform norm, the Weierstrass approximation theorem may be applied and $f$ can be approximated by polynomials at the desired precision.


 The topic  has been extensively studied by many authors, which were looking for a polynomial approximating  sequence having optimal Lebesgue constants in spaces of locally continuous functions equipped with uniform weighted norms \cite{mastrocco2001,mastrove,laurita,mastromilo2009,mastronota,mastrocco2010}.
Based on Laguerre zeros,  we mention the classical  Lagrange interpolating  polynomial  for which in \cite{mastrocco2001} it has been proved that
the associated weighted  Lebesgue constants diverge  as $n^{\frac 1 6}$. In the same paper a "remedium" was proposed  to this bad behavior, by introducing a suitable additional node. In this case,  sufficient and necessary conditions have been stated to get optimal Lebesgue constants growing  with  minimal order $\log n$.  A truncated version of the latter Lagrange polynomial has been introduced and studied in \cite{laurita,mastromilo2009} by involving only a very reduced number of function evaluations.  Also in this case we get  optimal  Lebesgue constants,  under suitable sufficient conditions.

%
%
%

In the present paper we introduce new approximation polynomials based on the same data of the truncated Lagrange-Laguerre interpolation but able to provide uniformly bounded Lebesgue constants and a near-best convergence order for any locally continuous function on the semiaxis. Such polynomials are defined by discretizing some VP means, namely delayed arithmetic means of the Fourier-Laguerre partial sums.

We recall that in literature such continuous VP means have been introduced as an useful tool to prove several theoretical results 
\cite{Mastro_Woula_ACTA,Mastro_Woula_JAT}.

In the Jacobi case it is known that  by discretizing VP means we get near-est approximation operators that preserve the polynomials up to certain degree, have uniformly bounded Lebesgue constants \cite{numeralg} and, in some cases, are interpolating too. Here we show that, in the Laguerre case, the interpolation property as well as the polynomial preservation are not ensured. But, all the same, we prove the uniform convergence at a near-best approximation rate with uniformly bounded Lebesgue constants, under suitable assumptions.

The outline of the paper is the following. In Section 2 the necessary notations and preliminary results are introduced. Section 3 concerns the definition and the properties of the new VP discrete approximation polynomials. Section 4 deals with the numerical experiments and, finally, Section 5 regards the conclusions.

\section{Notations and preliminary results}
\label{basic}

Along all the paper the notation $\C$ will be used several times to denote a positive constant having different values in different formulas. We will write $\C \neq \mathcal{C}(a,b,\ldots)$ in order to say that $\mathcal{C}$ is independent of the parameters $a,b,\ldots$, and $\mathcal{C} = \mathcal{C}(a,b,\ldots)$ to say that
$\mathcal{C}$ depends on $a,b,\ldots$. Moreover, if $A,B > 0$ are quantities depending on some parameters, we will write $A \sim B,$ if there exists an absolute constant $\C>0$, independent on such parameters,  such that
 $\C^{-1}B \leq A  \leq \C B$.

For all $\gamma\geq 0,$ we will consider the weight $\uu(x)=x^\gamma e^{-x/2},$ and the space $C_{\uu}$ of all locally continuous functions on the half line (i.e., continuous in $]0,\infty[$) such that
\[
 \lim_{x\to \infty}(f\uu)(x)=0 ,\qquad \mbox{and, in the case $\gamma>0$, also} \qquad
  \lim_{x\to 0}(f\uu)(x)=0.
\]
We will take the space $C_\uu$ equipped  with the norm
\begin{equation*}
    \|f\|_{C_{\uu}}= \|f\uu\|:=\sup_{x\geq 0}\left|(f\uu)(x) \right|,
\end{equation*}
and for any compact set $D\subset [0,+\infty)$,  we will use $\|f\|_D= \sup_{x\in D}|f(x)|$.

Note that functions in  $C_{\uu}$ may grow exponentially as $x\to +\infty$, and may have an algebraic singularity at the origin.

Denoting by  $\PP_n$ the space of the algebraic
polynomials of degree at most $n\ge 0$, we will consider the error of best approximation of $f\in C_{\uu}$ in $\PP_n$
\begin{equation}\label{best_approx}
 E_n(f)_{\uu}:=\inf_{P_n\in \PP_n}\| (f-P_n){\uu}\|,
\end{equation}
It is known that  \cite{mastroianni_ETNA,mastromilobook,MS1,MS2}
$$\lim_{n\to \infty} E_n(f)_{{\uu}}=0, \qquad \forall f\in C_\uu. $$

For \ smoother functions, we consider the Sobolev-type spaces of order $r\in \NN$
\[W_r(\uu)=\left\{ f\in C_{\uu} \ : \ f^{(r-1)} \in AC((0,+\infty)) \mbox{ and } \ \|f^{(r)}\varphi^r \uu\| < +\infty \right\},\]
where  $AC((0,+\infty))$ denotes the set of all functions which are  absolutely continuous on every closed subset of  $(0,+\infty)$ and $\varphi(x)=\sqrt{x}.$

We equip these spaces with the norm
\[\|f\|_{W_r(\uu)}:= \|f \uu\| + \|f^{(r)}\varphi^r \uu\|.\]

 For any function $f\in W_r(\uu)$ the following estimate holds \cite{viggiano,mastroianni_ETNA}
\begin{equation}\label{favard} E_n(f)_{{\uu}}\le \C \frac{\|f\|_{W_r(\uu)}}{\sqrt{n^r}},\quad \C\neq\C(n,f).
\end{equation}

\subsection{Truncated Lagrange interpolation at the zeros of Laguerre polynomials}
Let $\wal(x)=e^{-x}x^\alpha$ be the Laguerre weight of parameter $\alpha >-1$
and let  $\{p_n(\wal)\}_n$ be the
 corresponding sequence of orthonormal Laguerre polynomials with positive leading
coefficients.
Denoting by  $x_k:=x_{n,k},\ k=1,\dots,n,$ the zeros of $p_{n}(\wal)$ in increasing order, we recall that (see \cite{Sze})
\begin{equation}\label{xk}
\frac \C n< x_{1}<x_{2}<\ldots< x_{n}<4n+2\alpha-\C n^{\frac 13},\qquad \C\ne\C(n).
\end{equation}
In the sequel, we will be interested  in taking only a part of these zeros.
More precisely, throughout the paper,  the integer $j:=j(n)$ will denote the index of the  last zero of $p_n(\wal)$ we consider.  It is  defined by
\begin{equation} \label{nodomax}
j=\min_{k=1,2,..,n}\left\{k: x_{k}\ge 4n\rho \right\},\end{equation}
with $0<\rho<1$ arbitrarily fixed,

We recall that, inside the segment $[0,x_j]$, the distance between two consecutive zeros of $p_n(\wal)$ can be estimated as follows
\begin{equation*}
 \Delta x_k\sim \Delta x_{k-1} \sim \sqrt{\frac{x_{k}}{n}}, \qquad \quad \Delta x_{k}=x_{k+1}-x_{k},\qquad  k=1,2,\dots,j.
\end{equation*}

With $j$ defined in \eqref{nodomax}, let $L_{n+1}^*(\wal,f)$ be the truncated Lagrange polynomial \cite{mastromilo2009} defined by  \begin{equation}\label{lag_doppio}L_{n+1}^*(\wal,f,x):=\sum_{k=1}^j f(x_k)\ell_{n,k}(x),
\end{equation}
where
\begin{equation}\label{fund-L*}
\ell_{n,k}(x)=\frac{p_{n}(\wal,x)}{p'_{n}(x_{k})(x-x_{k})}\frac{4n-x}{4n-x_{k}},\qquad k=1,2,\dots, j,
\end{equation}
are the fundamental Lagrange polynomials  corresponding to the zeros of $p_n(\wal,x)(4n-x)$.

It is known that the norm of the operator  $L_{n+1}^*(\wal):f\in C_u\to L_{n+1}^*(\wal,f)\in C_u$, i.e. the weighted Lebesgue constant
\begin{equation}\label{norma_op_Lagrange}\|L_{n+1}^*(\wal)\|_{C_u}=\sup_{\|fu\|=1}  \|L_{n+1}^*(\wal,f)\uu\|=\sup_{x\ge 0}  \sum_{k=1}^j |\ell_{n,k}(x)| \frac{\uu(x)}{\uu(x_k)}, \end{equation}
diverges at least as $\log n$ when $n\to\infty$. To be more precise,  the following result holds \cite{mastromilo2009}

\begin{theorem}\cite{mastromilo2009} \label{th-C0_lag}
Let  $w_\alpha(x)=e^{-x}x^\alpha,\ \alpha>-1$ and $u(x)=e^{-\frac x 2}x^\gamma,\ \gamma\geq 0.$ Under the assumption
\begin{equation}\label{ipo_milov}
\max \left(0,\frac \alpha 2 +\frac 1 4 \right) \le  \gamma \le  \frac \alpha 2 + \frac 5 4,
\end{equation}
we have
\begin{equation}\label{leb_lag}
\| L_{n+1}^*(\wal,f)\uu\|\le \C \|f\uu\| \log n,
\quad \forall f\in C_{\uu},\qquad \C\ne\C(n,f).
\end{equation}
\end{theorem}
Moreover, as regards the approximation of $f\in C_u$ provided by the polynomial $ L_{n+1}^*(\wal,f)$ we recall the following

\begin{theorem}\cite[Th.2.2]{mastromilo2009}
For any $f\in C_u,$ under the assumption \eqref{ipo_milov}, we have
$$\|(f-L_{n+1}^*(\wal,f))u\|\le \C (\log n) \left[E_N(f)_u + e^{-An}\|fu\|\right],$$
where $N= \left\lfloor \left(\frac \rho {1+\rho}\right) n \right \rfloor\sim n$, and
$A$ and $\C$ are positive constants independent of $n$ and $f$.
\end{theorem}


\subsection{Ordinary and truncated Gauss-Laguerre rules}
According to the previous notation,  the   Gauss-Laguerre rule on the nodes in (\ref{xk}) is given by
\begin{equation}\label{gaussiana}
\int_0^{+\infty} f(x)\wal(x)dx=\sum_{k=1}^n f(x_{k} )\lambda_{n,k}(\wal)+R_{n}(f),\qquad n\in\NN,
\end{equation}
where $\lambda_{n,k}(\wal)=\left[\sum_{j=0}^np_j(x_k)\right]^{-1}$  are the Christoffel numbers w.r.t. $\wal,$ and the remainder term satisfies
\begin{equation}\label{order-G}
R_{n}(f)=0,\qquad \forall f\in \PP_{2n-1}.
\end{equation}
For any $f\in C_\uu$, under the assumption
$\alpha-\gamma>-1$, the following estimate holds
$$|R_{n}(f)| \le \ \C E_{2n-1}(f)_{\uu}, \qquad \C\neq \C(n,f).$$

In  \cite{mastromone} the ``truncated"  Gauss-Laguerre rule was introduced by using only the first $j$ zeros of $p_n(\wal)$, with $j$ defined in (\ref{nodomax}), i.e.
\begin{equation}\label{troncata}\int_0^{+\infty} f(x)\wal(x)dx=\sum_{k=1}^j f(x_{k} )\lambda_{n,k}(\wal)+R_{n,j}(f),\end{equation}
where
\begin{equation}\label{Rnj}
R_{n,j}(f)=0,\qquad \forall f\in \mathcal{P}_{n-1}^*:=\left\{p\in
\PP_{n-1}:p(x_{k})=0,\quad k>j\right\}.
\end{equation}

In \cite{mastromilo2009} the following estimate was proved under the assumption
$\alpha-\gamma>-1$
$$|R_{n,j}(f)| \le \C \ \left( E_N(f)_{\uu}+e^{-\mathcal{B}n}\|f \uu\|\right)
,\qquad \forall f\in C_\uu,$$
where $N= \left\lfloor \left(\frac \rho {1+\rho}\right) n \right \rfloor\sim n$,  and $\mathcal{B},\C>0$ are constants independent of $n,f$.

\subsection{Cesaro and de la Vall\'ee Poussin means}
For any given $n\in \NN_0=\{0,1,2,...\}$, let $S_n(\wal,f)\in \PP_n$ the $n-$th Fourier sum of $f$ in the orthonormal system $\{p_i(\wal)\}_{i=0,1,..}$, namely
$$S_n(\wal,f,x)=\sum_{i=0}^n c_i(f) p_i(\wal,x)$$
with Fourier coefficients defined as
\begin{equation}\label{ci}
c_i(f)=\int_0^{+\infty} f(t)p_i(\wal,t)\wal(t) dt.
\end{equation}
For all $n\in\NN=\{1,2,....\}$, we recall the Cesaro means are defined as arithmetic means of the first $n-1$ Fourier sums, i.e.
\begin{equation}\label{Cesaro}
\sigma_n(\wal,f,x)=\frac 1 {n} \sum_{k=0}^{n-1} S_k(\wal,f,x)=\sum_{i=0}^{n-1} \mu_{n,i}^m c_{i}(f) p_i(\wal,x),\end{equation}
In \cite{poiani} simple necessary and sufficient conditions have been stated in order that the map $\sigma_n(\wal):C_{\uu}\to C_{\uu}$ is uniformly bounded w.r.t. $n\in\NN$.
\begin{theorem}\label{th-Cesaro}
\cite[Theorems 1 and 5]{poiani}\newline
For all $\alpha>-1$, let  $w_\alpha(x)=e^{-x}x^\alpha$ and $u(x)=x^\beta\sqrt{w_\alpha(x)}$. The map $\sigma_n(\wal):C_{u}\to C_{u}$ is uniformly bounded w.r.t. $n\in\NN$, i.e.
we have
\begin{equation}\label{limitatezza}
\| \sigma_n(\wal,f)\uu\|\le \C \|f\uu\|,
\quad \forall f\in C_{\uu},\quad\forall n\in\NN,\qquad \C\ne\C(n,f),
\end{equation}
if and only if the following bounds are satisfied
\begin{equation}\label{ipo_Ces}
- \min \left\{\frac \alpha 2, \frac 1 4 \right\}< \beta <
1+\min\left\{\frac \alpha 2, \frac 14\right\},
\qquad\mbox{and}\qquad  - \frac 1 2\le \beta \le \frac 76
\end{equation}
\end{theorem}

Now, besides $n\in\NN$, let us take an additional parameter $m\in \NN$ such that $0<m<n$,  and, for a given function $f$, let us  consider the VP means, namely the following delayed arithmetic means of Fourier sums
\begin{equation}\label{VPcontinuo-mean}
\mathcal{V}_n^m(\wal,f,x)=\frac 1 {2m} \sum_{k=n-m}^{n+m-1} S_k(\wal,f,x),
\end{equation}
From a computational point of view, this is equivalent to
\begin{equation}\label{VPcontinuo}
\mathcal{V}_n^m(\wal,f,x)=\sum_{i=0}^{n+m-1} \mu_{n,i}^m c_{i}(f) p_i(\wal,x),
\end{equation}
that is a filtered version of $S_{n+m-1}(w_\alpha,f)$ with filter coefficients
$\mu_{n,i}^m$  defined as
\begin{equation}\label{muj}
\mu_{n,i}^m:=\left\{\begin{array}{cl}
1 & \mbox{if}\quad i=0,\ldots, n-m,\\ [.1in]
\displaystyle\frac{n+m-i}{2m} & \mbox{if}\quad
n-m+1\le  i\le  n+m-1.\\
\end{array}\right.
\end{equation}

The filtered VP polynomial $\mathcal{V}_n^m(\wal,f)\in \PP_{n+m-1}$ can be also rewritten in integral form as
\begin{equation}\label{VP_integrale}\mathcal{V}_n^m(\wal,f,x)=\int_0^{+\infty} f(t) v_{n}^m (x,t) w_\alpha(t) dt,\end{equation}
where
$$ v_{n}^m (x,t)=\frac 1 {2m}\sum_{r=n-m}^{n+m-1}K_r(x,t),\qquad\mbox{and}\qquad K_r(x,t)=\sum_{i=0}^r p_i(w_\al,x)p_i(w_\al,t)$$
are the VP and Darboux kernels, respectively.

Note that $v_n^m(x,t)=v_n^m(t,x)$ is a polynomial of degree $n+m-1$ in  both the variables that  can be represented as
$$v_n^m(x,t)=\sum_{i=0}^{n+m-1} \mu_{n,i}^m p_i(\wal,x)p_i(\wal,t).$$
An advantage of taking VP means instead of Cesaro means is given by the following invariance property
\begin{equation}\label{inva}\mathcal{V}_n^m(\wal,P)=P,\quad \forall P\in \PP_{n-m}.\end{equation}
Such preservation property allows us to investigate the approximation provided by $\mathcal{V}_n^m(\wal,f)$ in $C_u$, by estimating the norm of the operator map $\mathcal{V}_n^m(\wal):C_\uu\to C_\uu$, that is
\begin{equation}\label{norma_op_continuo}\|\mathcal{V}_n^m(\wal)\|_{C_u}=\sup_{\|fu\|=1}  \|\mathcal{V}_n^m(\wal,f)u\|=\sup_{x\ge 0}\left[ \uu(x) \int_0^{+\infty} |v_{n}^m (x,t)| \frac{w_\alpha(t)}{\uu(t)}dt\right]. \end{equation}

To this aim, we introduce the following notation concerning the degree parameters $m,n\in\NN$ defining the VP mean:
\begin{eqnarray}\label{nsim}
m\sim n & \mbox{iff}& m<n\le \C m,\qquad\qquad
\mbox{holds with $\C>1$ independent of $n,m$},\\
\label{napprox}
m\approx n & \mbox{iff}& m<c_1 m\le n\le c_2 m,\quad
\mbox{holds with $c_2\geq c_1> 1$ independent of $m,n$.}
\end{eqnarray}

Note that  $m\sim n$ agrees with the notation introduced in Section \ref{basic} while $m\approx n$ is a stronger condition that ensures $(n-m)\sim n$ holds too.

If we take the additional parameter $m$ such that $m\sim n$ then from Theorem \ref{th-Cesaro} we easily deduce the uniform boundedness of VP operator $\mathcal{V}_n^m(\wal)$ in $C_\uu$ for suitable choices of the weights $\wal$ and $\uu$. More precisely, we have the following
\begin{theorem}\label{th-C0}
Let the weight function  $w_\alpha(x)=e^{-x}x^\alpha$ and  $u(x)=e^{-\frac x 2}x^\gamma,$ satisfy the following bounds
\begin{equation}\label{ipo_new}
\max\left\{\frac\alpha 2 -\frac 14,\ 0\right\}<\gamma <\min\left\{\frac\alpha 2 +\frac 76,\ \alpha+1\right\}
\end{equation}
then, for any couple of positive integers $m \sim n$,  the map $\mathcal{V}_n^m(\wal):C_{\uu}\to C_{\uu}$ is uniformly bounded w.r.t. $n$ and $m$, i.e.
\begin{equation}\label{limitatezza}
\| \mathcal{V}_n^m(\wal,f)\uu\|\le \C \|f\uu\|,
\quad \forall f\in C_{\uu},\quad\forall m\sim n,\qquad \C\ne\C(n,m,f).
\end{equation}
\end{theorem}
{\it Proof of Theorem \ref{th-C0}}
Taking into account that
$$\mathcal{V}_n^m(\wal,f)=\frac{n+m}{2m}\sigma_{n+m}(w_\alpha,f)-\frac{n-m}{2m}\sigma_{n-m}(w_\alpha,f), $$
the statement follows by using the assumption  \eqref{nsim} that implies
$$\frac{n+m}{2m}\le \frac{\C+1}2,\qquad\mbox{and}\qquad  \frac{n-m}{2m}\le \frac{\C-1}2,
\qquad \C\ne\C(n,m),$$
and applying Theorem \ref{th-Cesaro} with $\beta=\gamma-\frac\alpha 2$, that satisfies (\ref{ipo_Ces}) by virtue of (\ref{ipo_new}).  \Proofend

\noindent Hence, concerning the error estimate  the following result can be deduced
\begin{theorem}\label{th-nearbest}
Let the weight function  $w_\alpha(x)=e^{-x}x^\alpha$ and  $u(x)=e^{-\frac x 2}x^\gamma,$ satisfy \eqref{ipo_new}. Then, for any $f\in C_u$ and for all $m\sim n$, we have
\begin{equation}\label{nearbest}
E_{n+m-1}(f)_{\uu}\le \|[f-\mathcal{V}_n^m(\wal,f)]\uu\|\le \C E_{n-m}(f)_{\uu},\qquad \C\ne \C(n,m,f),
\end{equation}
Moreover, if  $m \approx n$ then we get
\begin{equation}\label{lim-inf}
 \|f-\mathcal{V}_n^m(\wal,f)\|_{C_u}\sim E_n(f)_u, \qquad \forall f\in C_{\uu}
\end{equation}
with the constants in $\sim$ independent of $n,m,f$
\end{theorem}
{\bf Proof of Theorem \ref{th-nearbest}}

The first inequality in \eqref{nearbest} trivially follows from $\mathcal{V}_n^m(\wal,f)\in\PP_{n+m-1}$.

The second inequality in \eqref{nearbest} can be deduced by means of \eqref{inva} and \eqref{limitatezza} that yield
\[
\|[f-\mathcal{V}_n^m(\wal,f)]\uu\|\le \|(f-P)u\|+ \|\mathcal{V}_n^m(\wal, \ f-P)u\|\le\C \|(f-P)u\|, \qquad \forall P\in\PP_{n-m}
\]
and the statement follows by taking, at both the sides, the infimum w.r.t. $P\in\PP_{n-m}$.
\Proofend

\section{VP filtered discrete approximation}
In the previous section, we have seen that for any $f\in C_u$ and suitable choices of $\wal$, the polynomials $\{\mathcal{V}_n^m(\wal,f)\}_{m\approx n}$ provide a near--best approximation of $f$. However, their construction  requires the knowledge of the Fourier coefficients of the function $f$. To overcome this problem, in this section
we introduce the polynomial ${V}_n^m(\wal,f)$, that we define as a discrete version of $\mathcal{V}_n^m(\wal,f)$.
To be more precise, with $j=j(n)$ defined in \eqref{nodomax}, we apply the  Gauss-Laguerre rule in \eqref{troncata} to the coefficients $c_i(f)$ in \eqref{ci}, getting
\begin{equation}\label{cij}
c_{n,i}^{(j)}(f)=\sum_{k=1}^j f(x_k)p_i(\wal,x_k)\lambda_{n,k}(\wal)
\end{equation}
Hence, similarly to \eqref{VPcontinuo} we define the following discrete VP filtered approximation  polynomials

\begin{equation}\label{VPdiscreto}
{V}_n^m(\wal,f,x)=\sum_{i=0}^{n+m-1} \mu_{n,i}^m c_{n,i}^{(j)}(f) p_i(\wal,x),\qquad 0<m<n
\end{equation}
that can be equivalently written as
\begin{equation}\label{VPdiscreto-phi}
{V}_n^m(\wal,f,x)=\sum_{k=1}^j f(x_k)\Phi_{n,k}^m(x),\qquad x\ge 0
\end{equation}
with $j=j(n)$ defined in \eqref{nodomax}, and
\begin{equation}\label{phi}
\Phi_{n,k}^m(x)=\lambda_{n,k}(\wal)v_n^m(x,x_k)=
\lambda_{n,k}(\wal)\sum_{i=0}^{n+m-1} \mu_{n,i}^m p_i(\wal,x_k)p_i(\wal,x).
\end{equation}

\begin{figure}[!h]
\centering
\includegraphics [height=5.5cm,width=9.5cm] {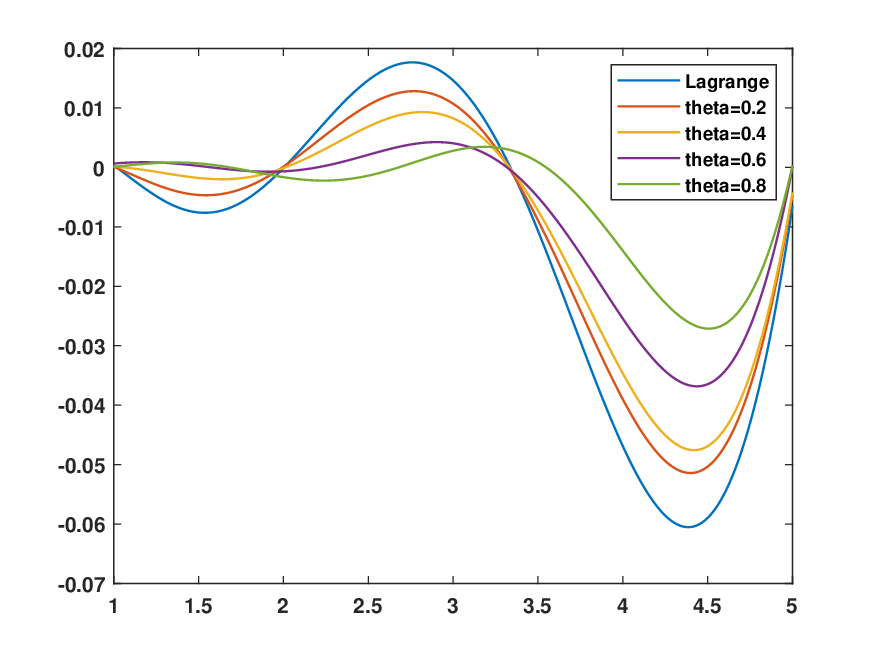}
\caption{\label{fund1} Plots of the fundamental polynomials $\ell_{n,k}(x)$ and $\Phi_{n,k}^m(x)$ for $\alpha=-0.5$, $n=15$, $k=7$, and $m=\lfloor \theta n\rfloor$}
\end{figure}

In Fig. \ref{fund1} we show, as an example,  the graphics of the
fundamental polynomial   $\Phi_{n,k}^m(x)$, for $n$ and $k$ fixed, and  $m$ varying. As we can see for increasing $m$
the oscillations of $\Phi_{n,k}^m(x)$ progressively dampen. The major oscillating behavior is presented  by the fundamental Lagrange polynomial $\ell_{n,k}(x)$ that are plotted in the same figure.

We remark that, similarly to $\mathcal{V}_n^m(\wal,f)$, for any function $f$ we have ${V}_n^m(\wal,f)\in \PP_{n+m-1}$, but in the case $f=P\in\PP_{n-m}$ the invariance property (\ref{inva}) does not apply to $V_n^m(\wal, f)$.

It can be easily proved that the norm of the operator map ${V}_n^m(\wal):C_\uu\to C_\uu$ takes the following form
\begin{equation}\label{norma_op_discreto}
\|{V}_n^m(\wal)\|_{C_u}=\sup_{\|fu\|=1}  \|{V}_n^m(\wal,f)\|=\sup_{x\ge 0} \sum_{k=1}^j \left|\Phi_{n,k}^m(x)\right|\frac{\uu(x)}{\uu(x_k)}. \end{equation}
As the number of nodes increases, the behavior of such norm is an important measure of the conditioning, due to the linearity of $V_n^m(\wal)$.
Under the same assumptions of Theorem \ref{th-C0}, the next theorem states that also the discrete filtered VP operator ${V}_n^m(\wal):C_{\uu}\to C_{\uu}$ is uniformly bounded  w.r.t any $n\sim m$.

\begin{theorem}\label{th-C0-dis}
If the weights $w_\alpha(x)=e^{-x}x^\alpha$ and $u(x)=e^{-\frac x 2}x^\gamma$ satisfy (\ref{ipo_new}) then, for any  $m \sim n$,  the map ${V}_n^m(\wal):C_{\uu}\to C_{\uu}$ is uniformly bounded w.r.t. $n$ and $m$, i.e.
\begin{equation}\label{limitatezza-dis}
\| {V}_n^m(\wal,f)\uu\|\le \C \|f\uu\|,
\quad \forall f\in C_{\uu},\qquad \C\ne\C(n,m,f).
\end{equation}
\end{theorem}
{\bf Proof of Theorem \ref{th-C0-dis}}

By Remez-type polynomial inequality \cite{mhaskar}, since ${V}_n^m(\wal,f)\in \PP_{n+m-1}$, we have
\begin{eqnarray*}\| {V}_n^m(\wal,f)\uu\|&\le & \C \max_{x\in \mathcal{A}_n} \left |{V}_n^m(\wal,f,x)\uu(x)\right|,\\ \quad \mathcal{A}_n&=&[a_{n,m},b_{n,m}],\quad a_{n,m}\sim \frac c {n+m-1},\quad  b_{n,m}\sim 4(n+m)\end{eqnarray*} 
by which
\begin{equation}\label{bound}\| {V}_n^m(\wal,f)\uu\|\le \C \|f\uu\| \left(\max_{x\in \mathcal{A}_n}\sum_{k=1}^j \left| \Phi_{n,k}^m(x)\right| \frac{\uu(x)}{\uu(x_k)}\right)\end{equation}
Hence, our aim is to prove that
\begin{equation}\label{tesi}
\Sigma_{n,m}(x):=\sum_{k=1}^j \left| \Phi_{n,k}^m(x)\right| \frac{\uu(x)}{\uu(x_k)}\le \C
,\qquad \forall x\in \mathcal{A}_n,\qquad \C\ne\C(n,m,x)
\end{equation}
Start from \eqref{phi} and use \begin{equation}\label{stimalambda}\lambda_{n,k}(\wal)\sim \Delta x_k \wal (x_k),\quad \Delta x_k=x_{k+1}-x_k,\end{equation} to obtain
\begin{equation}\label{stimaSigma}
\Sigma_{n,m}(x)\le \uu(x)\sum_{k=1}^j\Delta x_k \frac{\wal (x_k)}{\uu(x_k)} \left|v_n^m(x,x_k)\right|, \qquad \forall x\in \mathcal{A}_n,
\end{equation}
with $j$ defined in \eqref{nodomax}.

Now, let us recall that the following Marcinkiewicz-type inequality was proved in \cite[Lemma 2.5]{mastrove} for an arbitrary polynomial $P\in \PP_{nl}$ with $l\in\NN$
\begin{equation}\label{Marci}
\sum_{k=1}^j \Delta x_k |(Pu)(x_k)| \le \C \int_{x_1}^{4m\rho_1} |P(y)| u(y) dy, \quad \C\neq \C(m,P).
\end{equation}
where $0<\rho<\rho_1<1$, and $\C \neq \C(n)$.

By applying \eqref{Marci} to the the polynomial $P(y)=v_n^m(x,y)\in \PP_{2n}$ and the weights in \eqref{stimaSigma},  we get
\begin{equation}\label{sigma1}
\Sigma_{n,m}(x)\le \C \int_{x_1}^{4n\rho_1} \frac{\wal (t)}{\uu(t)} \left|v_n^m(x,t)\right| dt\le \C \int_{0}^{+\infty} \frac{\wal (t)}{\uu(t)} \left|v_n^m(x,t)\right| dt.\end{equation}
Consequently, by  \eqref{norma_op_continuo} and Theorem \ref{th-C0}, we deduce
$$\Sigma_{n,m}(x)\le \C \uu(x)\int_{0}^{+\infty} \frac{\wal (t)}{\uu(t)} \left|v_n^m(x,t)\right| dt\le \C \|\mathcal{V}_n^m(\wal)\|_{C_\uu}\le \C ,$$
i.e. (\ref{tesi}) holds and the theorem follows.\Proofend

Finally, by the next theorem we state that for suitable weight functions, the VP polynomials $V_n^m(\wal,f)$, similarly to their continuous version $\mathcal{V}_n^m(\wal,f)$ uniformly converge to $f$ in the space $C_u$, with a near--best convergence order. Hence, in conclusion, the VP filtered polynomials $V_n^m(\wal,f)$ are simpler to compute than the polynomials $\mathcal{V}_n^m(\wal,f)$ and share with the latter a near--best convergence rate in the spaces $C_u$.
\begin{theorem}\label{th-convergenza}
Let the weight function  $w_\alpha(x)=e^{-x}x^\alpha$ and  $u(x)=e^{-\frac x 2}x^\gamma,$ satisfy \eqref{ipo_new}. Then, for any $f\in C_u$ and for all $m\sim n$, we have
\begin{equation}\label{err-dis}
\| (f-{V}_n^m(\wal,f))\uu\|\le \C \left( E_{q}(f)_{\uu}+ e^{-An}\|f\uu\| \right),\quad q=\min\left\{n-m,\ \left\lfloor n\frac{\rho}{(1+\rho)}\right\rfloor\right\},
\end{equation}
with $0<\rho<1$ given by (\ref{nodomax}) and $\C,A>0$ independent of $n,m,f.$
Moreover, if  $m \approx n$ then we get
\begin{equation}\label{lim-inf1}
 \|f-V_n^m(\wal,f)\|_{C_u}\sim E_n(f)_u, \qquad \forall f\in C_{\uu}
\end{equation}
with the constants in $\sim$ independent of $n,m,f$.
\end{theorem}
{\bf Proof of Theorem \ref{th-convergenza}}

Let $P^*\in \PP_{n-m}$ be the polynomial of best approximation of $f\in C_\uu$, i.e., let
\begin{equation}\label{P*}
E_{n-m}(f)_u=\inf_{P\in \PP_{n-m}}\|f-P\|_{C_u}=\|f-P^*\|_{C_u}.
\end{equation}
Hence,  consider the following decomposition
$$\| (f-{V}_n^m(\wal,f))\uu\|\le \| [f-P^*]\uu\|+ \| {V}_n^m(\wal,f-P^*)\uu\| +
\|[P^*- {V}_n^m(\wal,P^*)]\uu\|.$$
By applying Theorem \ref{th-C0-dis} and \eqref{P*} we deduce
\begin{equation}\label{stima}
\| (f-{V}_n^m(\wal,f))\uu\|\le
\C E_{n-m}(f)_{\uu} + \|[P^*- {V}_n^m(\wal,P^*)]\uu\|\end{equation}
In order to estimate the last addendum, we note that using \eqref{inva} and \eqref{gaussiana}--\eqref{order-G} we can write
\begin{eqnarray*}
P^*(x)&=& \mathcal{V}_n^m(\wal,P^*,x)=\int_0^{+\infty}P^*(t)v_n^m(x,t)\wal(t)dt\\
&=& \sum_{k=1}^n P^*(x_k)v_n^m(x,t)\lambda_{n,k}(\wal)= \sum_{k=1}^n P^*(x_k)\Phi_{n,k}^m(x)
\end{eqnarray*}
On the other hand, by \eqref{VPdiscreto-phi}, we have
\[
V_n^m(\wal, P^*,x)=\sum_{k=1}^j P^*(x_k)\Phi_{n,k}^m(x), \qquad j=\min\{k \ :\ x_k\ge 4n\rho\}
\]
with $0<\rho<1$ arbitrarily fixed. Thus, for any $x\ge0$, it is
\begin{eqnarray}
\nonumber
&& 
|P^*(x)-V_n^m(\wal,P^*,x)|\uu(x)=\left|\sum_{k=j+1}^n P^*(x_k)\Phi_{n,k}^m(x)\right|\uu(x)\\
\nonumber
&\le&\sum_{k=j+1}^n |P^*(x_k)| |\Phi_{n,k}^m(x)|\uu(x)\\
\nonumber
&\le & \max_{k> j}|(P^*\uu)(x_k)|\left( \sum_{k=j+1}^n \left|\Phi_{n,k}^m(x)\right|\frac{\uu(x)}{\uu(x_k)}\right)\\
\nonumber
&\le & \|P^*\uu\|_{(4n\rho,+\infty)}\left(\max_{x\ge 0} \sum_{k=1}^n \left|\Phi_{n,k}^m(x)\right|\frac{\uu(x)}{\uu(x_k)}\right)\\
\nonumber
&=& \|P^*\uu\|_{(4n\rho,+\infty)}\ \|{V}_n^m(\wal)\|_{C_u}\\
\label{stima-P*}
&\le& \C \|P^*\uu\|_{(4n\rho,+\infty)}
\end{eqnarray}
where in the last two lines we have applied \eqref{norma_op_discreto} and Theorem \ref{th-C0-dis}, respectively.

Finally, we recall that for sufficiently large $n$, we have for any fixed positive $0<\rho<1$ \cite[Proposition 2.1]{mastromilo2009}
\begin{equation}\label{eq-inf-fin}
\|f\uu\|_{[4n\rho,+\infty)}\le E_N(f)_\uu + \C e^{-A n}\|f\uu\|,\qquad \forall f\in C_\uu
\end{equation}
where $N=\left\lfloor \frac{\rho }{1+\rho}n\right\rfloor,$  and $\C,A$ positive constants independent of $m,f.$

Consequently, by \eqref{P*} and \eqref{eq-inf-fin}, we get
\begin{eqnarray*}
\|P^*\uu\|_{(4n\rho,+\infty)}&\le& \|(P^*-f)\uu\|_{(4n\rho,+\infty)}+\|f\uu\|_{(4n\rho,+\infty)}\\
&\le& \C E_{n-m}(f)_u +  E_M(f)_\uu + \C e^{-A n}\|f\uu\|
\end{eqnarray*}
and the statement follows combining last bound with \eqref{stima} and \eqref{stima-P*}.\Proofend

\section{Numerical experiments}
In this section, we propose some tests to assess the performance of the truncated VP approximation by making comparisons with the truncated Lagrange approximation, under different points of view. To be more precise we will propose tests about
\begin{enumerate}
\item
the approximation of functions with different smoothness degree,
\item
the behavior of the weighted Lebesgue constants for different choices of the parameters $\alpha$ and $\gamma$.
\end{enumerate}

\noindent In what follows we will denote by  $Y$ a sufficiently large uniform mesh of points in a finite interval  $(0,a)$ with  $a>0$ sufficiently large. For any $x\in Y$, we consider the following weighted and unweighted error functions

$$e_{n,m}^{VP}(f,x)=|V_n^m(\wal, f, x)-f(x)|\uu(x),\quad \widetilde e_{n,m}^{VP}(f,x)=|V_n^m(\wal, f, x)-f(x)|$$
$$e_{N}^{Lag}(f,x)=|L_{N+1}^*(\wal,f,x)-f(x)|\uu(x),\quad \widetilde e_{N}^{Lag}(f,x)=|L_{N+1}^*(\wal,f,x)-f(x)|$$
where $L_{N+1}^*(\wal,f,x)$ is the truncated interpolating polynomial interpolating $f$ at the nodes $\{x_{N,k}\}_{k=1}^N \cup \{4N\}$, being $\{x_{N,k}\}_{k=1}^N$ the zeros of $p_N(\wal)$. We will consider  two  cases: $N=n$ and $N=n+m-1$. In the first case, $L_{n+1}^*(\wal,f)$ and $V_n^m(\wal, f)$ require  the same number of function evaluations, but they have different degrees. In the second case,  $L_{n+m}^*(\wal,f)$ and $V_n^m(\wal, f)$ have  the same degree, but  $L_{n+m}^*(\wal,f)$ uses a greater number of function evaluations.

Moreover, we will compute the maximum weighted errors given by
\begin{equation}\label{Err-test}
\mathcal{E}_{n,m}^{VP}:=\max_{y\in Y}\biggm(e_{n,m}^{VP}(f,y)\biggm),\qquad\mbox{and}
\mathcal{E}_{N+1}^{Lag}:=
 \max_{x\in Y}\biggm(e_{N}^{Lag}(f,y)\biggm),
\end{equation}
where $N\in\{n,n+m-1\}$

\subsection{Approximation of functions}
In the following we will consider six examples concerning the approximation of six test functions having different smoothness. \newline
For increasing values of $n$, we report the values of the maximum weighted errors in tables where the first two columns display the values of  $n$ and $m$  defining the VP polynomial $V_n^m(\wal, f)$, the 3rd column shows  the number of function evaluations for computing both the polynomials $V_n^m(\wal, f)$ and  $L_{n+1}^*(\wal,f)$, in the 4th and 5th columns  we find
 the  errors $\mathcal{E}_{n,m}^{VP}$ and $\mathcal{E}_{n+1}^{Lag}$, and, finally, in the 6th and 7th columns
the number of function evaluations required for  $L_{n+m}^*(\wal,f)$  and the corresponding error $\mathcal{E}_{n+m}^{Lag}$ are reported.\newline
Moreover, in order to show the quality of the pointwise approximation, we plot the unweighted error functions $\tilde e^{VP}_{n,m}(f,x)$ and $\tilde e^{Lag}_{n+1}(f,x)$ for fixed $n$ and $m$, letting x variable in a meaningful compact interval.
\begin{es}\label{es1} Let
\[\alpha=-0.4,\quad \gamma=0.05\quad
f_1(x)=\exp{(x/4)}\]
 The function $f_1$ is analytic and, according to the theory, the maxima weighted errors quickly decrease, as shown in Table \ref{tab:b4}. The unweighted error functions $\tilde e^{VP}_{n,m}(f_1,x)$ and $\tilde e^{Lag}_{n+1}(f_1,x)$ are plotted in Figure \ref{figex5}.
\begin{table}[ht]
\begin{center}
\begin{tabular}{|l|l|c|l|l|c|l|} \hline
$n$  & $m$  & \# f.eval & $\mathcal{E}_{n,m}^{VP}$  & $\mathcal{E}_{n}^{Lag}$ & \# f.eval &    $\mathcal{E}_{n+m}^{Lag}$                \\ \hline
 20 	 &    6 	 &   19 	 & 8.00e-08 	 & 1.80e-05    &   25 	  & 2.21e-07 	      \\ \hline
 70 	 &    7 	 &   58 	 & 2.70e-14 	 & 2.76e-14    &   61 	  & 8.90e-14 	 \\ \hline
 120 	 &   12 	 &   79 	 & 5.83e-14 	 & 5.77e-14    &   83 	  & 9.14e-14 	 \\ \hline
 170 	 &   51 	 &   96 	 & 4.49e-14 	 & 5.86e-14    &  110 	  & 3.02e-13 	 \\ \hline
 220 	 &   22 	 &  110 	 & 6.25e-14 	 & 6.55e-14    &  115 	  & 7.07e-13 	 \\ \hline
 \end{tabular}
\caption{\label{tab:b4} Maxima weighted errors induced by Lagrange and VP polynomials for Example \ref{es1}
}
\end{center}
\end{table}
 \begin{figure}[!h]
\centering
{\includegraphics[height=5.5cm,width=8.5cm]{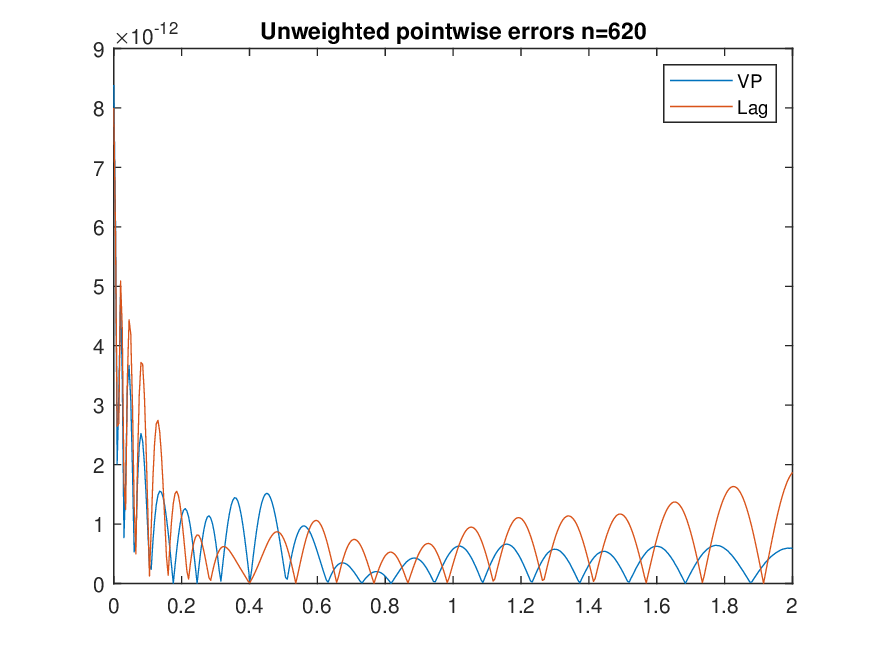}}
{\caption{\label{figex5} Example \ref{es1}: Plots of $\tilde e^{VP}_{n,m}(f_1,x)$ and $\tilde e^{Lag}_{n+1}(f_1,x)$ for $n=620$, $m=186$. }}
\end{figure}   	
\end{es}  	
\begin{es}\label{es2}
Consider
\begin{eqnarray*}
\alpha=0.5,\quad \gamma=0.5\quad
f_2(x)=\frac{|x-1|^{\frac{11}2}}{100+x^2}\end{eqnarray*}
In this case, $f_2$ belongs to an intermediate space between $W_1(u)$ and $W_2(u)$.   According to \cite[(5.2.7)]{mastromilobook},  the theoretical convergence order is $\mathcal{O}\left(\frac 1 {(n-m)^\frac{11}4}\right)$ for $\mathcal{E}_{n,m}^{VP}$, and
$\mathcal{O}\left(\frac{\log N}{N^\frac{11}4}\right)$ for $\mathcal{E}_{N}^{Lag}$.
As shown in Table \ref{tab:b2}, all the numerical errors  agree with the theoretical estimates and are comparable as a consequence of the very slow divergence of the logarithmic factor. In particular, VP approximation  provides a slightly better approximation than the interpolating polynomial in the case $N=n$. This trend  is little bit reversed when $N=n+m$.

However, in Figure \ref{figex2}, the graphics of the pointwise error function taken for $n=420,\ m=126$, shows that  $V_n^m(\wal,f_2)$ provides a pointwise approximation better than the one offered by $L_{n+1}^*(\wal,f_2)$.
\begin{table}[!ht]
\begin{center}
\begin{tabular}{|l|l|c|l|l|c|l|} \hline
$n$  & $m$  & \# f.eval & $\mathcal{E}_{n,m}^{VP}$  & $\mathcal{E}_{n}^{Lag}$ & \# f.eval &    $\mathcal{E}_{n+m}^{Lag}$                \\ \hline
  20 	 &    6 	 &   19 	& 1.10e-04 	      	  & 4.17e-02          &   25      &    2.92e-02                   \\ \hline
 220 	 &   22 	 &   96 	 & 9.78e-08 	    & 9.23e-08      &  101 	  &  	      5.17e-08                      \\ \hline
 420 	 &   42 	 &  134 	 & 1.27e-08 	    & 1.38e-08      &  140 	  & 	       8.54e-09                     \\ \hline
 620 	 &  310 	 &  163 	 & 1.54e-09 	    & 2.36e-09      &  200 	  & 	       1.42e-09                    \\ \hline
 820 	 &   82 	 &  188 	 & 2.47e-09 	    & 2.53e-09      &  197 	  &         1.63e-09                    \\ \hline
1020 	 &  102 	 &  210 	 & 1.31e-09 	    & 1.31e-09      &  220 	  &         1.01e-09                      \\ \hline
1220 	 &  122 	 &  230 	 & 7.28e-10 	    & 7.39e-10      &  241 	  &         5.75e-10                      \\ \hline
\end{tabular}
\caption{\label{tab:b2} Maxima weighted errors induced by interpolating  and VP polynomials for Example \ref{es2}.}
\end{center}
\end{table}

\begin{figure}[!h]
\centering
\includegraphics [height=5.5cm,width=7.5cm]{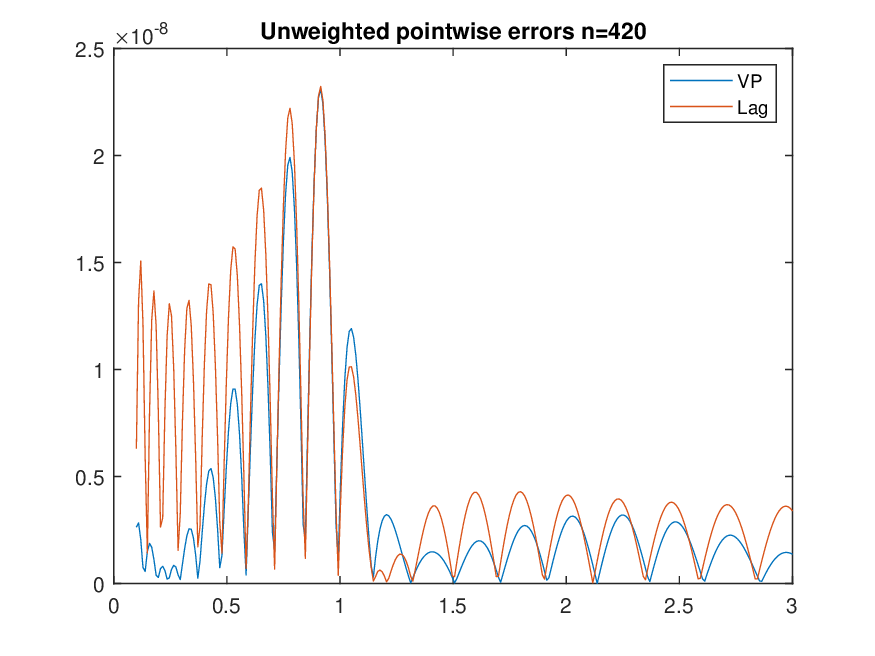}
\caption{\label{figex2} Example \ref{es2}: Plots of $\tilde e^{VP}_{n,m}(f_2,x)$ and $\tilde e^{Lag}_{n+1}(f_2,x)$ for $n=420$, $m=126$.}
\end{figure}
\end{es}  	

\begin{es}\label{es3}
Consider the case
\begin{eqnarray*}
&& f_3(x)=\frac 1{1+100(x-3)^2},\quad \alpha=-0.4,\quad \gamma=0.05\quad\end{eqnarray*}
\begin{table}[ht]
\begin{center}
\begin{tabular}{|l|l|c|l|l|c|l|} \hline
$n$  & $m$  & \# f.eval & $\mathcal{E}_{n,m}^{VP} $  & $\mathcal{E}_{n}^{Lag}$ & \# f.eval &    $\mathcal{E}_{n+m}^{Lag}$                \\ \hline
    20 	     &    6 	 &   17 	 & 2.24e-01 	 	  & 2.25e-01 	  &   20 	  & 1.13e-01                                                            \\ \hline
    220 	 &   22 	 &   64 	 & 9.04e-02 	   & 8.95e-02 	 &   67 	  & 1.57e-01                       	    	                                 \\ \hline
    420 	 &  378 	 &   89 	 & 6.98e-02 	   & 8.34e-02 	 &  123 	  & 8.31e-02                       		                                 \\ \hline
    620 	 &  558 	 &  108 	 & 6.21e-02 	   & 7.05e-02 	 &  149 	  & 4.12e-02                       		                                \\ \hline
    820 	 &  738 	 &  124 	 & 3.72e-02 	   & 4.93e-02 	 &  172 	  & 2.56e-02                       	                                 \\ \hline
   1020 	 &  102 	 &  139 	 & 5.86e-02 	   & 5.82e-02 	 &  146 	  & 5.66e-02                      	                                    \\ \hline
   1220 	 & 1098 	 &  152 	 & 2.84e-02 	   & 3.64e-02 	 &  209 	  & 2.87e-02                                                           \\ \hline
\end{tabular}
\caption{\label{tab:b3} Maxima weighted errors induced by interpolating  and VP polynomials for Example \ref{es3}.}
\end{center}
\end{table}

Note that $f_3\in W_{r}(\uu)$, $\forall r>0$, and the slow convergence we see in Table \ref{tab:b3} depends on the growth of the norms in the error estimates. For instance, with $r=5$, $\|f_3\|_{W_5(\uu)}\sim 0.5e+12$. Similarly to the previous example, the maxima errors attained by  VP and Lagrange approximation are almost comparable among them. However, also in this case in Figure \ref{figex3} are displayed the absolute pointwise errors taken for $n=420,\ m=378$, and except a small range around the peak point  $x=3$, the errors by   $V_n^m(\wal,f)$ are smaller than those by $L_{n+1}^*(\wal,f)$.

\begin{figure}[!h]
\centering
\includegraphics[height=5.5cm,width=8.5cm]{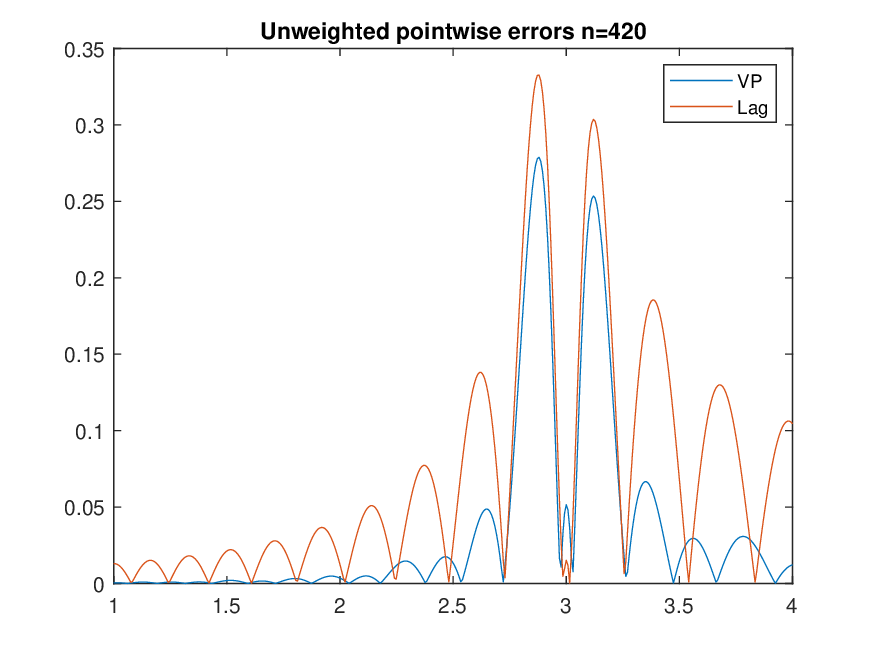}
\caption{\label{figex3} Example \ref{es3}: Plots of $\tilde e^{VP}_{n,m}(f_3,x)$ and $\tilde e^{Lag}_{n+1}(f_3,x)$ for $n=420$ and $m=378$ }
\end{figure}
\end{es}

\begin{es}\label{es4} Let be
\begin{eqnarray*}
\alpha=0.6,\quad \gamma=0.6\quad
f_4(x)=|\cos(\pi-x)|\end{eqnarray*}
\begin{table}[ht]
\centering
\begin{tabular}{|l|l|c|l|l|c|l|} \hline
$n$      & $m$       & \# f.eval & $\mathcal{E}_{n,m}^{VP}$ & $\mathcal{E}_{n}^{Lag}$ & \# f.eval &    $\mathcal{E}_{n+m}^{Lag}$                \\ \hline
  20 	 &    2 	 &   19 	 & 1.37e-01 	  	  & 1.36e-01 	    &   21 	  & 1.16e-01         \\ \hline
 220 	 &   22 	 &   79 	 & 3.24e-02 	  	  & 3.25e-02 	 &   83 	  & 3.95e-02 	          \\ \hline
 420 	 &   84 	 &  108 	 & 2.57e-02 	  	  & 2.63e-02 	 &  121 	  & 1.65e-02 	      \\ \hline
 620 	 &   62 	 &  137 	 & 1.75e-02 	  	  & 1.73e-02 	 &  141 	  & 1.23e-02 	     \\ \hline
 820 	 &   82 	 &  155 	 & 1.69e-02 	  	  & 1.67e-02 	 &  152 	  & 1.63e-02 	     \\ \hline
1020 	 &  102 	 &  173 	 & 1.76e-02 	  	  & 1.74e-02 	&  182 	  & 1.25e-02 	      \\ \hline
1220 	 &  488 	 &  190 	 & 8.31e-03 	  	  & 8.01e-03 	&  225 	  & 9.35e-03 	       \\ \hline
\end{tabular}
\caption{\label{tab:b10} Maxima weighted errors induced by interpolating  and VP polynomials for Example \ref{es4}.}
\end{table}
In this case we have $f_4\in W_{1}(\uu)$,  and the degree of approximation is very poor, since, according to \eqref{favard}, VP and Lagrange approximation rates are $\mathcal{O}\left(\frac 1 {\sqrt{n-m}}\right), \mathcal{O}\left(\frac{\log n}{\sqrt{n}}\right)$, respectively. The numerical results displayed in Table \ref{tab:b10} and Figure \ref{figex4} agree with the theory and, as the previous cases, we see that the maxima errors are almost comparable but, concerning the local approximation, VP approximation error behaves better away from the singularity $x=\pi/2$.
 \begin{figure}[!h]
\centering
{\includegraphics[height=5.5cm,width=8.5cm]{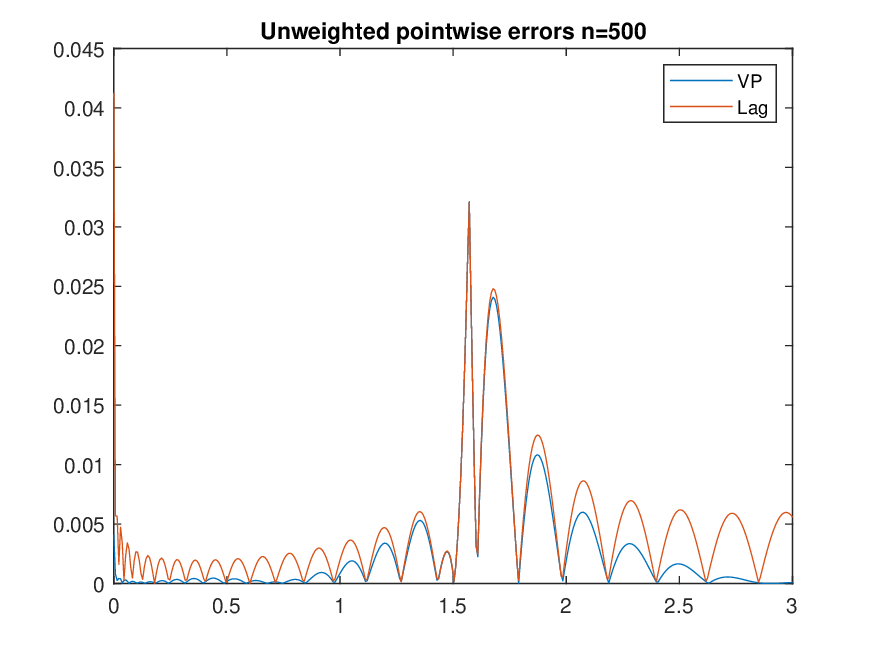}}
\caption{\label{figex4} Example \ref{es4}: Plots of $\tilde e^{VP}_{n,m}(f_3,x)$ and $\tilde e^{Lag}_{n+1}(f_3,x)$ for $n=500$ and $m=150$.}
\end{figure}   	
\end{es}  	
\begin{es}
\label{ex_gibbs}
Now consider the following function plotted in Figure  \ref{fun} (left--hand side)
$$f(x)=\frac 1 {1+100(x-0.5)^2}+\frac 1{1+1000\sqrt{x^2+0.5}},\quad   \gamma=0.5$$
To evidence the improvement in the pointwise approximation achieved by VP polynomials, on the right--hand side of Fig. \ref{fun} we report the plots of both Lagrange and VP polynomials when $n=400$, $\alpha=0$ and $m=200.$
\begin{figure}[!h]
\begin{center}
{\includegraphics [height=5.5cm,width=9.5cm] {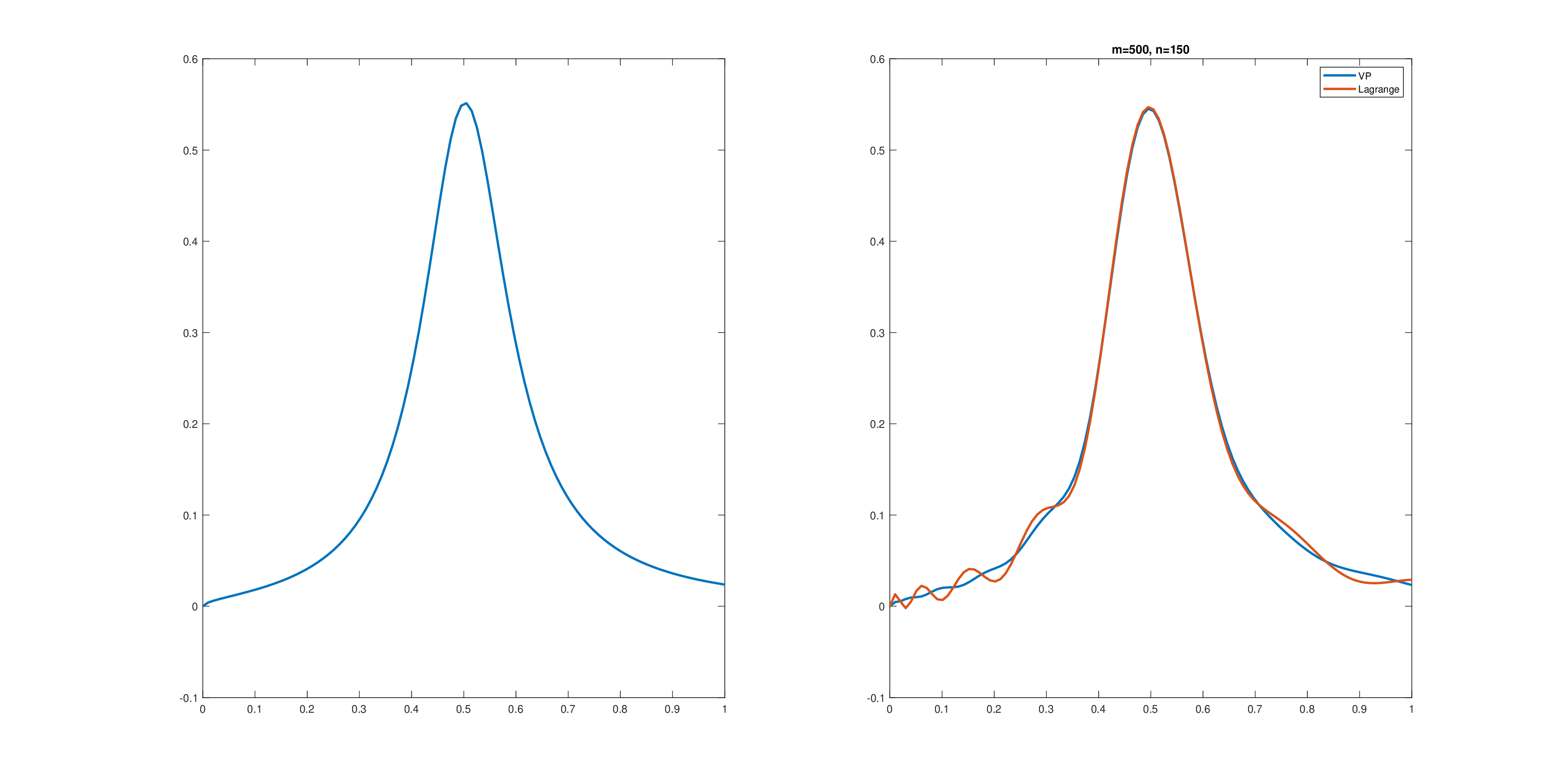}}
\hspace{5mm}
{\caption{\label{fun} Example 5: graphic of the function $f_5$ (left) and plots of the VP and Lagrange polynomials (right) for $\alpha=0$, $n=400$ and $m=200.$}}
\end{center}
\end{figure}
\end{es}
\begin{es}\label{es_bv}
This test highlights the feature offered by VP polynomials when we want to  approximate functions of bounded variation having a jump discontinuity. Consider the function
$$f_6(x)=\begin{cases} x, & x<1\\ x+2, &x\ge 1\end{cases}$$

Figure \ref{bv_function} displays the graphics of $f_6 \uu$ with $\uu(x)=e^{-\frac x 2 }\sqrt{x}$ (black line) and the approximations by the weighted Lagrange polynomial
$\uu L_{n+1}^*(\wal,f_6)$ (orange line) and the weighted $\uu V_n^m(\wal, f_6)$ (blue line),  for $n=620$ and $m=310$.

As expected, the Gibbs phenomenon appears. In fact, as the graphics in Figure \ref{bv_function} show, the Lagrange polynomial presents  damped oscillations that persist even when moving away from the jump while the  VP polynomial induces a reduction of the Gibbs phenomenon as $m$ increases. This is confirmed by Figure \ref{errore_bv} displaying the unweighted absolute error induced by Lagrange polynomial  and VP polynomial
\begin{figure}[!h]
\centering
{\includegraphics [height=6.5cm,width=14.5cm] {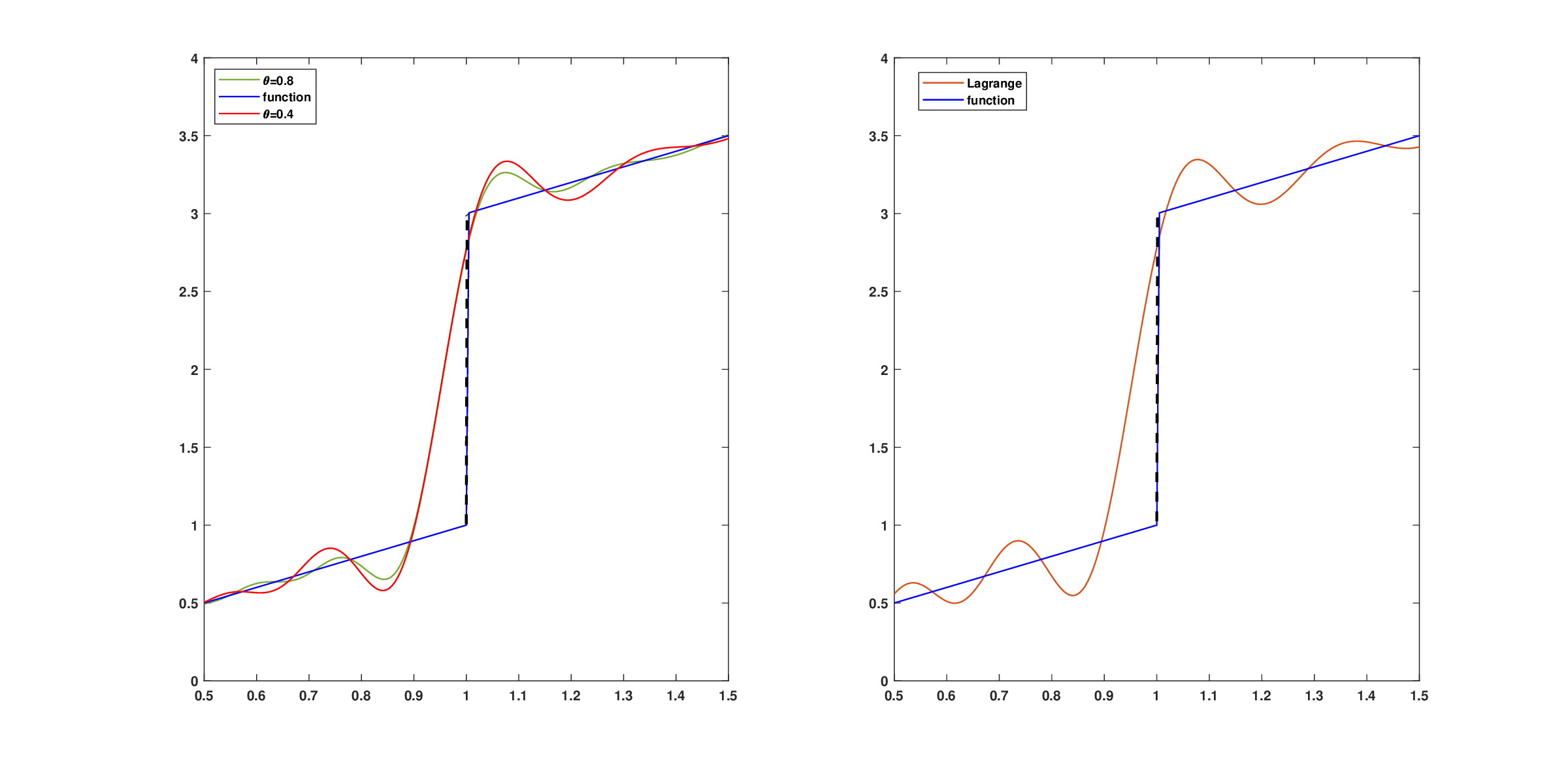}}
\hspace{5mm}
{\caption{\label{bv_function} $n=620$, $\alpha=0.5$, $\gamma=0.5$: the function $f_6\uu$ and weighted VP polynomials for $m=298$,  $m=496$ (left), and the weighted Lagrange polynomial of degree $n$ (right). }}
\end{figure}
\begin{figure}[!h]
\centering
{\includegraphics [height=7.5cm,width=12.5cm] {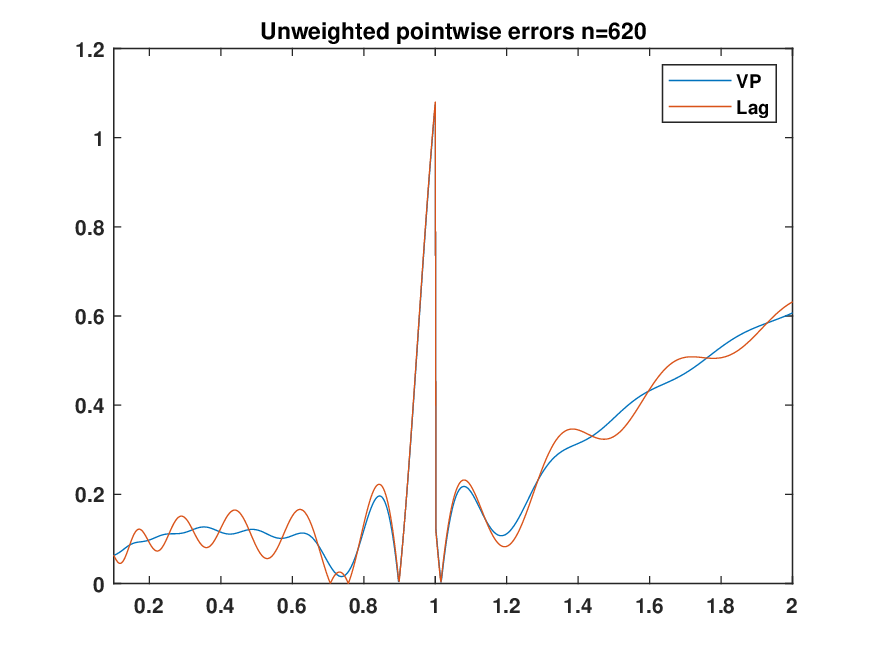}}
\hspace{5mm}
{\caption{\label{errore_bv} $n=620$, $m=310$,  $\alpha=0.5$. }}
\end{figure}
\end{es}

\subsection{Behavior  of the weighted Lebesgue constants}
\begin{figure}[!h]
{\includegraphics [height=5.5cm,width=7.5cm] {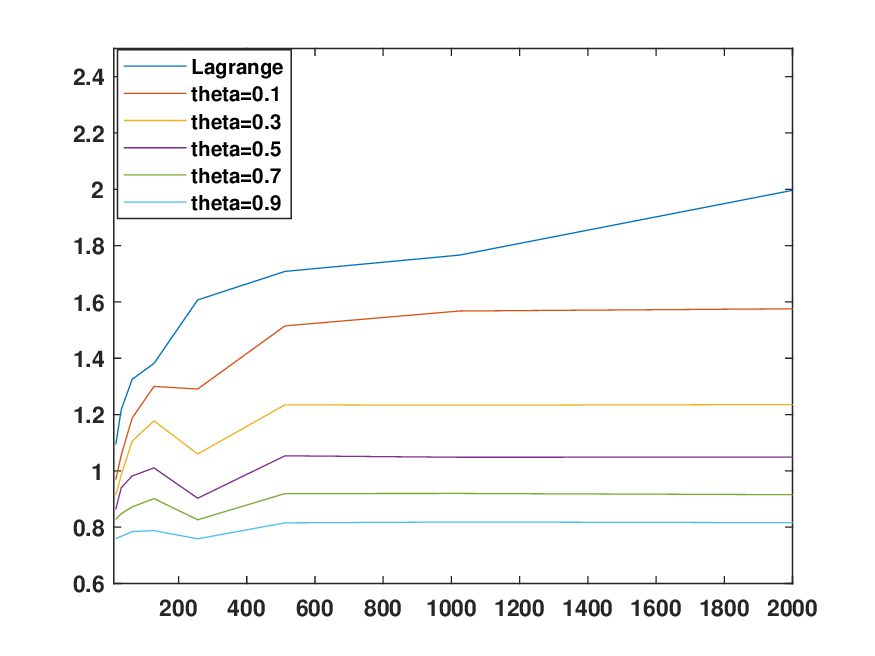}\hspace{5mm} \includegraphics [height=5.5cm,width=7.5cm] {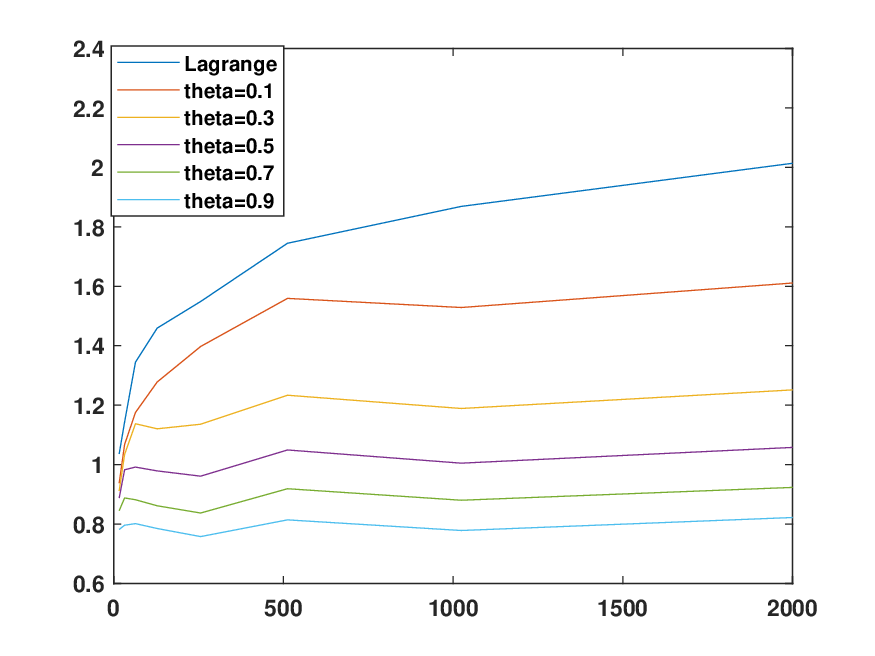}}
\hspace{5mm}
{\caption{\label{leb1}Lebesgue constants for $\alpha=\gamma=0.5$ (left), for for $\alpha=\gamma=1$ (right). }}
\end{figure}

\begin{figure}[!h]
{\includegraphics [height=5.5cm,width=7.5cm] {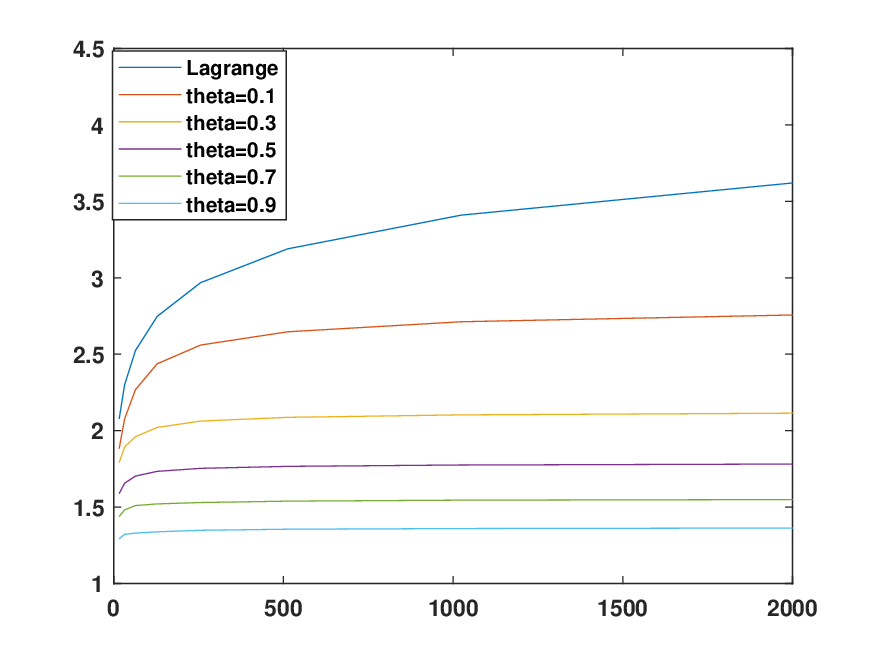}\hspace{5mm}\includegraphics [height=5.5cm,width=7.5cm] {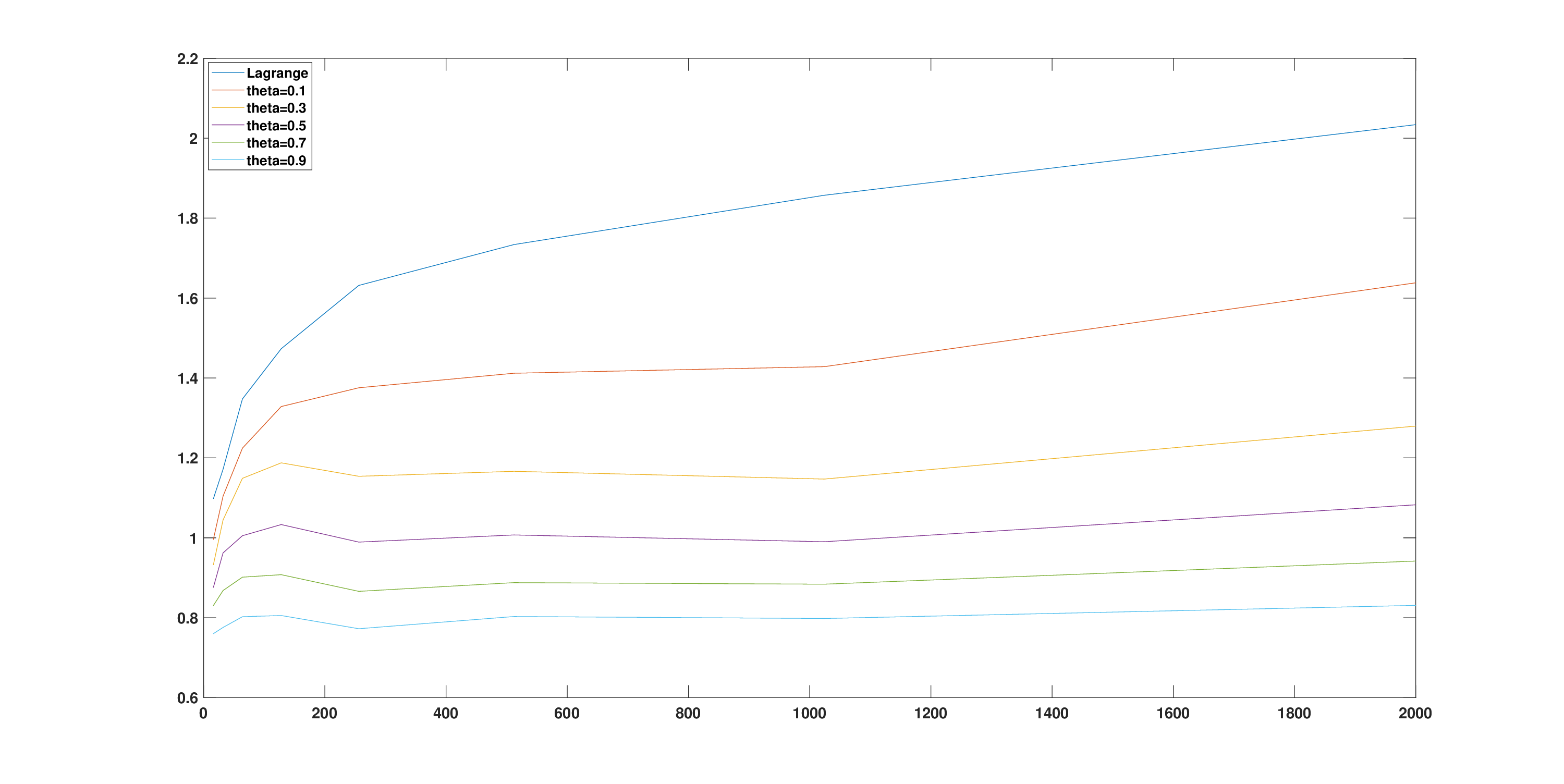}}
{\caption{\label{leb3}Lebesgue constants for $\alpha=-0.5, \gamma=0$ (left), $\alpha=0.7, \gamma=0.6$ (right). }}
\end{figure}

In what follows we show  the behavior of the
Lebesgue constants related to Lagrange and VP approximation for different choices of $\alpha$ and $\gamma$, as $n$ increases. From the theory sufficient conditions relating $\alpha$ and $\gamma$ are known to ensure the optimal logarithmic grow in the Lagrange case (cf. \eqref{ipo_milov}),  and the uniform boundedness in the VP case (cf. \eqref{ipo_new}). In the latter case, it is also required that $m\approx n$ to get a near best approximation order. To this aim, we take $m=\lfloor n\theta \rfloor$, being
 $\theta\in (0,1)$ chosen as   $\theta\in \{0.1, 0.2,\dots,0.9.\}$. Hence for these values of $\theta$ and $n$ increasing, the behavior of the Lebesgue constants is displayed
in Figures \ref{leb1}-\ref{leb3} for several couples of $(\alpha,\gamma)$ both in the ranges to obtain optimal  Lebesgue constants for VP and Lagrange operators.

\begin{figure}[!h]
{\includegraphics [height=5.5cm,width=7.5cm] {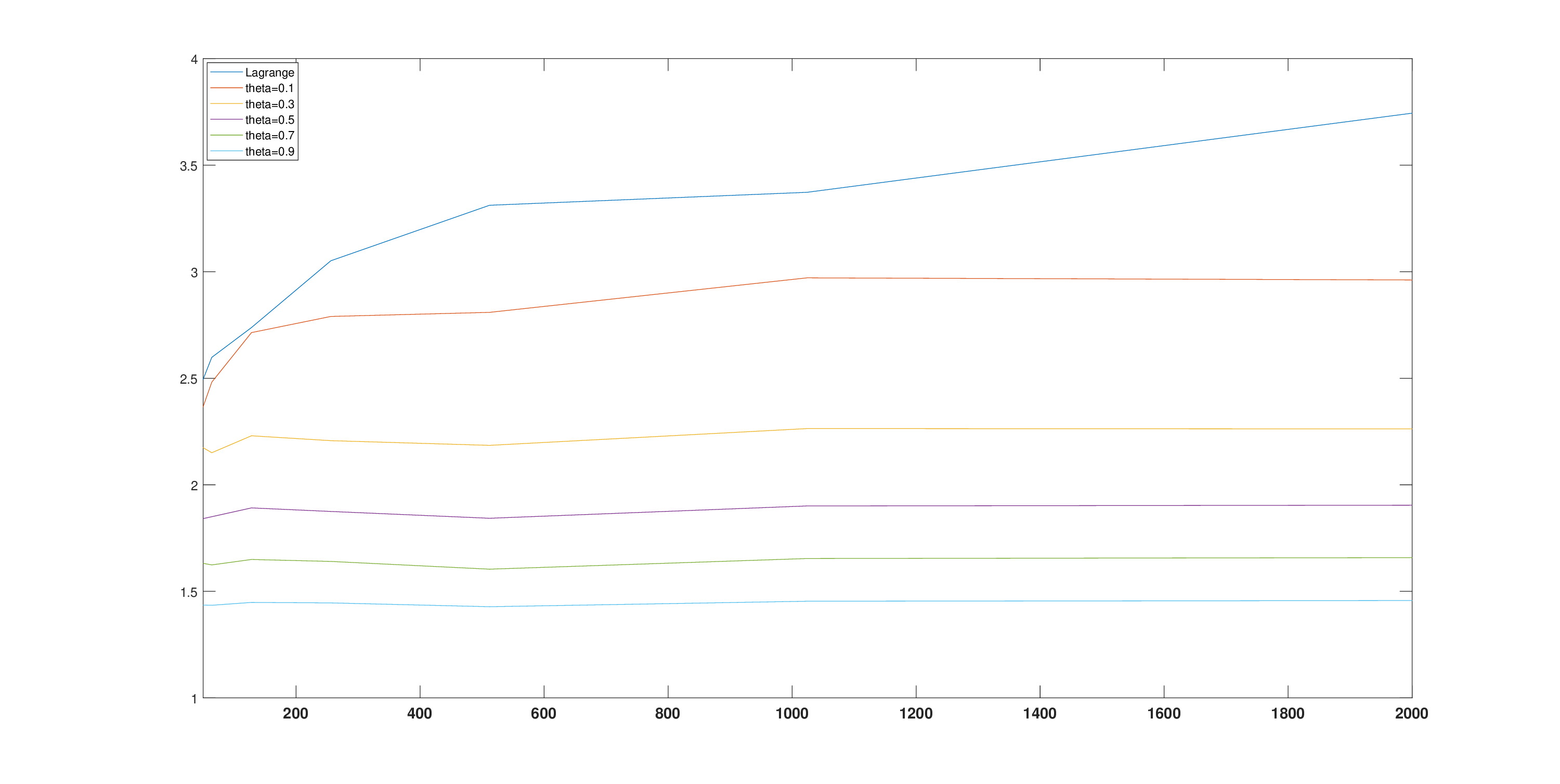} \hspace{5mm} \includegraphics [height=5.5cm,width=7.5cm] {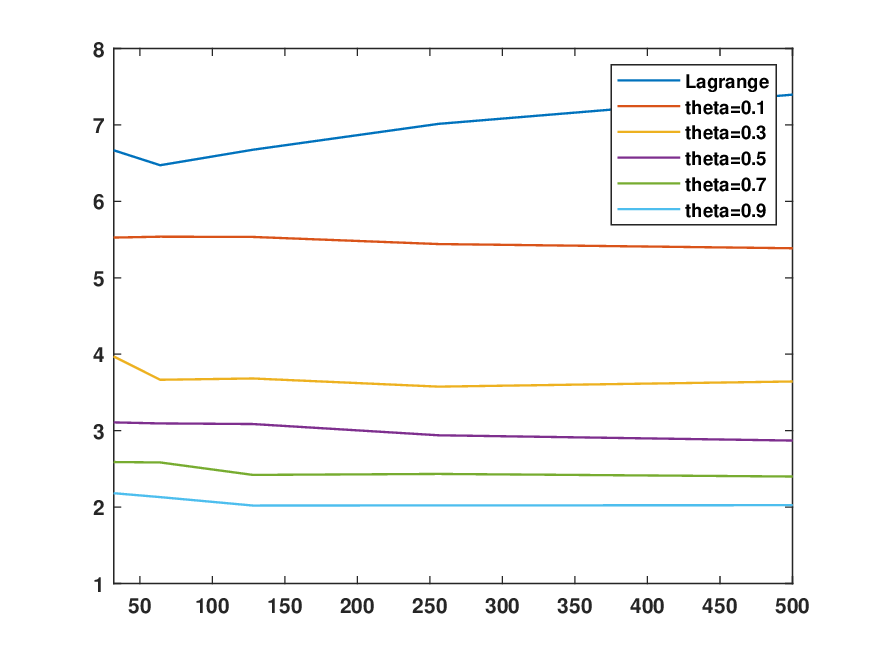}}
\hspace{5mm}
{\caption{\label{leb6}Lebesgue constants for $\alpha=0.5, \gamma=1.43$ (left), $\alpha=0, \gamma=5/4$ (right). }}
\end{figure}

In conclusion, in Figure \ref{leb6} we consider the case when $\alpha$ and $\gamma$ do not satisfy the bounds in \eqref{ipo_new} concerning VP approximation, but they satisfy those in \eqref{ipo_milov} concerning Lagrange.
In this case, even if not predicted by Theorem \ref{th-C0-dis} we see the Lebesgue constants of VP operator are still bounded for any $\theta$.

Empirically we have detected that  the same behavior occurs for other possible choices of couples of $(\alpha,\gamma)$ s.t. $\frac \alpha 2 +\frac 7 6<\gamma\le  \frac\alpha 2 +\frac 5 4.$
This observation leads us to conjecture that the hypotheses of Theorem  \ref{th-C0-dis} can be extended by replacing the upper bound  $\min\{\frac\alpha 2 +\frac 76, \alpha +1\}$ by $\frac\alpha 2 +\frac 5 4$.

\section{Conclusions}

We have introduced and studied  a new truncated polynomial operator  $V_n^m(\wal):C_u\to C_u$, involving the values of $f$ at only the first $j<<n$ zeros of the  Laguerre polynomial $p_n(w_\alpha)$. Differently from $L_{n+1}^*(\wal)$, we have proved that  the map $V_n^m(\wal):C_u\to C_u$ is uniformly bounded under suitable conditions  on the weights $\wal$ and $u$ and in view of this theoretical improvement,  the maxima errors attained by VP  approximation are better on average than those achieved by Lagrange one. However,   the differences between the numerical errors in norm are not so large, since  the factor $\log n$ does not have a strong impact in the error estimate.    The situation is instead quite different about  the pointwise errors, especially 
   in approximating  smooth functions presenting some  peaks or  isolated pathologies (see Examples \ref{es4}-\ref{ex_gibbs}). Indeed, in such  cases the additional parameter $m$ can be suitably modulated to reduce the Gibbs phenomenon. We observe a drastic reduction of the  oscillations in VP polynomials  as $m$ increases, and even more w.r.t the Lagrange polynomial, the latter swinging more than anyone. Such reduction is valuable also in Example \ref{es_bv} where $f_6$ is a bounded variation function. In these cases it is well–known that Lagrange polynomials present overshoots and oscillations not only close to the singularities, but also in the smooth part of the function. The graphics in Fig. \ref{bv_function} show how the Gibbs phenomenon can be  reduced by VP polynomials, for suitable values of the parameter $\theta\in (0, 1).$
We conclude observing that the   range where $\alpha,\gamma$ vary for the boundedness of VP polynomials, is partially  overlapping with that the norm of the Lagrange operator  logarithmically  behaves (compare \eqref{ipo_milov} and \eqref{ipo_new}), and in any case the  VP range is larger than the Lagrange one. However,  numerical results (see Fig. \ref{leb6}) allow us to conjecture that the range determined by Thm. \ref{th-C0-dis} about  the parameters $\alpha,\gamma$  could actually be stretched. More precisely, we conjecture that the hypothesis  $\max\{0,\frac \alpha 2-\frac 14\}<\gamma\le \frac \alpha 2+\frac 54$ could replace \eqref{ipo_new}.

\vspace{0.5cm}

\noindent\textbf{Acknowledgments.}
This work was partially supported by INdAM-GNCS.
The research has been accomplished within RITA (Research ITalian network on Approximation)

\bibliographystyle{plain}
\bibliography{biblio}

\end{document}